\documentclass{amsart}
\usepackage{epsfig}
\usepackage{psfrag}

\newtheorem{thm}{Theorem}[section]

\newtheorem{clm}[thm]{Claim}

\theoremstyle{definition}

\theoremstyle{remark}

\numberwithin{equation}{section}

\begin{document}

\title{A note on Kneser-Haken finiteness} 

\author{David Bachman}
\address{Mathematics Department, California Polytechnic State University}
\email{dbachman@calpoly.edu}


\date{September 24, 2002 and, in revised form, October 14, 2002.}


\begin{abstract}
Kneser-Haken Finiteness asserts that for each compact 3-manifold $M$ there is an integer $c(M)$ such that any collection of $k>c(M)$ closed, essential, 2-sided surfaces in $M$ must contain parallel elements. We show here that if $M$ is closed then twice the number of tetrahedra in a (pseudo)-triangulation of $M$ suffices for $c(M)$. 
\end{abstract}
\maketitle


\section{Introduction}

A crucial step in Kneser's proof of the existence of prime factorizations \cite{kneser:29} of 3-manifolds is the following Theorem:

\begin{thm}
\label{kneser}
Suppose that $M$ is a compact 3-manifold. Then there is an integer $c(M)$ such that if $\mathcal S=\{S_1,...,S_k\}$ is any collection of $k > c(M)$ pairwise disjoint 2-spheres in $\mbox{Int}(M)$ then the closure of some component of $M-\mathcal S$ is a punctured 3-cell.
\end{thm}

In \cite{haken:61.2} W. Haken generalized this Theorem to collections of incompressible surfaces. The precise result is the following:

\begin{thm}
\label{haken}
Suppose that $M$ is a compact, irreducible 3-manifold. Then there is an integer $c(M)$ such that if $\mathcal S=\{S_1,...,S_k\}$ is any collection of $k > c(M)$ pairwise disjoint, closed, 2-sided incompressible surfaces in $\mbox{Int}(M)$ then the closure of some component of $M-\mathcal S$ is a product.
\end{thm}

Theorems \ref{kneser} and \ref{haken} have become standard topics of introductory courses on the topology of 3-manifolds. Together they are often referred to as ``Kneser-Haken Finiteness". Several proofs exist, each giving different quantities for $c(M)$. Most use some function of the number of tetrahedra $t$ in a triangulation of $M$. For example, in \cite{haken:68} Haken uses $61t$ for $c(M)$. In \cite{hempel}, one of the standard texts used in introductory 3-manifolds courses, J. Hempel uses $c(M)=6t+\mbox{dim}H_1(M)+\mbox{dim}H_1(M,\mathbb Z_2)$. In lecture notes available on the world wide web \cite{hatcher:01} A. Hatcher uses the bound $c(M)=8t+\mbox{dim}H_1(M,\mathbb Z_2)$. A. Casson also has a proof in a set of widely distributed, handwritten lecture notes from a course in China. His bound is similar: $6t+2\mbox{dim}H_2(M,\mathbb Z_2)$. W. Neumann took up the task of typing Casson's notes \cite{neumann:01}, but gives a different proof for Kneser's theorem. He credits this proof to Dave Bayer, who finds that if $M$ is closed the quantity $6t$ is sufficient. Finally, in \cite{jr:98} Jaco and Rubinstein mention, without proof, that in the case that $M$ is closed one can take $5t$ for $c(M)$. 

In this note we show that if $M$ is closed then we can take $c(M)$ to be just $2t$. To prove this we build on the ``classical" proof by borrowing techniques from Bayer's proof and the proof of Jaco and Rubinstein mentioned in the previous paragraph, as well as introducing some new ideas of our own. 

Of course, for Kneser and Haken the important part of Theorems \ref{kneser} and \ref{haken} is the existence of {\it some} bound. The bound itself was unimportant. However, with the aid of computers people have begun tabulating 3-manifolds with relatively few tetrahedra and studying their properties. For such manifolds knowledge of the {\it best} bound for Kneser-Haken finiteness becomes relevant.  

In the interest of brevity, and by virtue of the fact that Kneser-Haken finiteness is such a standard topic of introductory courses, we will assume here that the reader is familiar with at least one proof of Theorem \ref{kneser} which uses the theory of normal surfaces in triangulated 3-manifolds. 

The author would like to thank his thesis advisor, Cameron McA. Gordon, for his initial exposure to Kneser-Haken Finiteness and comments on the first draft of this paper. The author would also like to thank Saul Schleimer for pointing him towards Bayer's proof, and further comments on an earlier draft. Finally, the author is greatly indebted to William Jaco for explaining his proof of Kneser-Haken finitness, providing Claim \ref{twisted}, and making several other helpful suggestions. 

\section{Theorem and Proof.}

In the classical setting of Theorems \ref{kneser} and \ref{haken} $M$ is assumed to be triangulated. Recently {\it pseudo}-triangulations have gained in popularity, so we assume this is the setting for our main result. Note that a pseudo-triangulation is one in which simplices are allowed to meet each other (and themselves) in a collection of lower dimensional simplices. This allows, for example, a single 3-simplex to be glued to itself along two of its faces. 

\begin{thm}
Let $M$ be a closed, irreducible, orientable 3-manifold equipped with a (pseudo)-triangulation with $t$ 3-simplices. If $\mbox{Int}(M)$ contains a collection of $k$ closed, non-parallel, essential, 2-sided surfaces then 
\[k\le2t\]
\end{thm}

\begin{proof}
Let $\mathcal S=\{S_1,...,S_k\}$ be a collection of essential, non-parallel, 2-sided surfaces in a compact 3-manifold $M$ with triangulation $T$. Let $T^i$ denote the $i$-skeleton of $T$. Assume that $\mathcal S$ has been isotoped to have least weight ({\it i.e.} $|\mathcal S \cap T^1|$ is assumed to be minimal) and is normal with respect to $T$. 

Let $M^*$ denote $M$ with a regular neighborhood of $T^0$ removed. Note that $\partial M^*$ is a collection of spheres, each made up of normal triangles. If two surfaces are parallel in $M^*$ then they are parallel in $M$, so the elements of $\mathcal S$ are non-parallel in $M^*$. If $X$ is any subset of $M$ then let $X^*=X \cap M^*$. 

Let $\Delta$ be a 3-simplex of $T$. A component of $(\partial \Delta)^* -\mathcal S$ is called {\it good} if it is an annulus. Otherwise it is {\it bad}. Call a component of $\Delta ^* -\mathcal S$ a {\it good cell} or a {\it bad cell} according as its intersection with $(\partial \Delta) ^*$ is good or bad. Note that there are precisely two possibilities for bad cells. These are a {\it truncated tetrahedron} ({\it i.e.} a subset of a tetrahedron bounded by 4 non-parallel normal triangles) or a {\it truncated prism} (as outlined in bold in Figure \ref{f:badquad}). Good cells are of the form $(\mbox{triangle}) \times I$ or $(\mbox{quadrilateral}) \times I$. 

Cut $M^*$ along $\mathcal S$ and call a component of $M^*-\mathcal S$ {\it good} if it is made up of good cells and {\it bad} otherwise. Each good component is an $I$-bundle over some surface. Note that no such $I$-bundle is trivial as that would imply that two elements of $\mathcal S$ were parallel or that one was parallel to the link of a vertex of $T^0$. 

A surface on the boundary of a component of $M^*-\mathcal S$ which is not the link of a vertex will be referred to as a {\it remnant} of $\mathcal S$. As there are $k$ elements of $\mathcal S$, and each such surface is 2-sided, there must be $2k$ remnants. A remnant is good or bad according as it appears on the boundary of a good or bad component of $M^*-\mathcal S$. Let $g$ denote the number of good remnants and $b$ the number of bad ones, so that $2k = g+b$.

For each good remnant there is a twisted $I$-bundle embedded in $M$. Each such $I$-bundle contributes one $\mathbb Z_2$ direct summand to $H_1(M;\mathbb Z_2)$. The proof of the following claim was related to the author by William Jaco. 

\begin{clm}
\label{twisted}
The rank of $H_1(M;Z_2)$ is bounded above by $t + 1$.
\end{clm}

\begin{proof}
Let $e$ and $v$ denote the number of edges and vertices of $T$. As $M$ is closed we have $e = t + v$. Let $\Gamma$ be a maximal tree of $T^1$, $e_\Gamma$ the number of edges in $\Gamma$, and $e_N$ the number of edges not in $\Gamma$. Then $1 = \chi(\Gamma) = v-e_\Gamma$. Combining this with the equalities $e = t + v$ and $e_\Gamma +e_N=e$ we obtain $e_N = t + 1$. The proof is complete by noting that the rank of $H_1(M;Z_2)$ is bounded above by $e_N$. 
\end{proof}


We say a triangle or quadrilateral (or, in general, a normal disk) is {\it bad} if it lies in the intersection of a bad remnant with a bad cell. Let $s$ and $q$ denote the number of bad triangles and bad quadrilaterals appearing on the boundaries of the bad cells of $M^* -\mathcal S$. Note that every tetrahedron contains at most four bad triangles and at most two bad quadrilaterals, so that $s \le 4t$ and $q \le 2t$. 

\begin{clm}
\label{two}
Every bad remnant contains at least two bad normal disks.
\end{clm}

\begin{proof}
Every bad remnant meets at least one bad cell. Hence it contains at least one bad triangle or one bad quadrilateral. Suppose first that there is a bad remnant, $R$, which contains a bad quadrilateral $Q$ and no other bad normal disk. Note that $Q$ lies on the boundary of a truncated prism. Such a cell has a pair of hexagonal faces. Since $Q$ is the unique bad normal disk of $R$, and $M$ is orientable, we conclude that these hexagonal faces are identified as in Figure \ref{f:badquad}. Inspection of this figure also shows that there is a disk $D$ with boundary on $R$. If $\partial D$ is essential on $R$ then we contradict the assumption that $R$ was incompressible. On the other hand, if $\partial D$ is inessential on $R$ then we contradict the assumption that $\mathcal S$ has minimal weight, since compressing $R$ along $D$ produces an isotopic normal surface which meets $T^1$ in fewer points. 

        \begin{figure}[htbp]
        \psfrag{Q}{$Q$}
        \psfrag{D}{$D$}
        \vspace{0 in}
        \begin{center}
        \epsfxsize=3 in
        \epsfbox{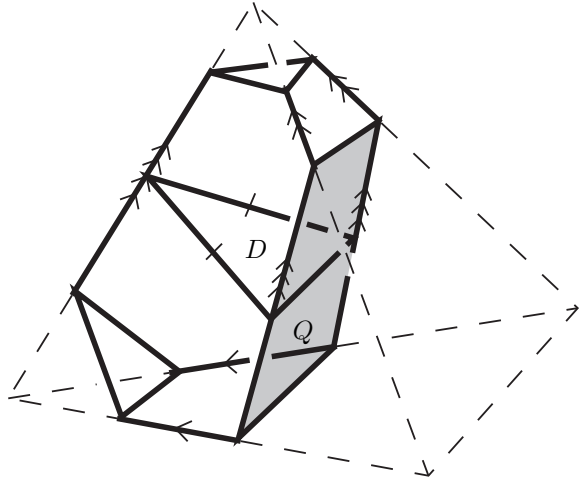}
        \caption{A single bad quad, $Q$, in a truncated prism.}
        \label{f:badquad}
        \end{center}
        \end{figure}

If $R$ contains a bad triangle $\Delta$ then there are two cases. Suppose first that $\Delta$ appears on the boundary of a truncated tetrahedon. Such a bad cell has four hexagonal faces where it meets other bad cells. As all three sides of $\Delta$ lie on such a face, and at most two can be identified, we conclude that $\Delta$ is glued to at least one other bad normal disk. 

The final case is when $\Delta$ lies on the boundary of a truncated prism. But if $R$ contains no other bad normal disks then again it must be the case that the hexagonal faces of the prism are glued as in Figure \ref{f:badquad}. Once again we conclude that either there is a compression for some other remnant (the one that contains the quadrilateral of the same prism) or the collection $\mathcal S$ did not have minimal weight.
 
\end{proof}

It follows immediately from Claim \ref{two} that the number of bad remnants $b$ is bounded above by $\frac{s+q}{2}$. Combining the above inequalities we now obtain:

\begin{eqnarray*}
2k &=& g+b\\
&\le& t+1+\frac{s+q}{2}\\
&\le & t+1+\frac{4t+2t}{2}\\
&=& 4t+1
\end{eqnarray*}
We conclude that $k \le 2t+\frac{1}{2}$. However, since $k$ is an integer it must be the case that $k \le 2t$
\end{proof}

\bibliographystyle{alpha}
\bibliography{finite}

\end{document}